\documentclass[10pt]{amsart}

\usepackage{amscd, amssymb}

\usepackage{graphicx}
\usepackage{psfrag}

\newenvironment{custom}[1]{\medskip {\bf #1:  } \em  }{\\}
\newtheorem{thm}{Theorem}[section]

\newtheorem{lem}[thm]{Lemma}
\newtheorem{cor}[thm]{Corollary}

\def\length #1{||#1  ||}
\newcommand{\ignore}[1]{\index{ignore}}
\newcommand{\st}{{\pi_{\mu}}}
\newcommand{\R}{{\mathbb R}}
\newcommand{\rp}{{\mathbb R}^+}

\newcommand{\Z}{\mathbb{Z}}
\newcommand{\N}{\mathbb{N}}

\renewcommand{\O}{{\mathcal O}}
\newcommand{\ibar}{\bar{\imath}}

\newcommand{\bpf}{\vspace{+10pt} \noindent {\bf{Proof.}} }
\newcommand{\epf}{\qed \vspace{+10pt}}

\renewcommand{\a}{\alpha}
\renewcommand{\b}{\beta}
\renewcommand{\d}{\delta}
\renewcommand{\th}{\theta}
\newcommand{\Th}{\Theta}
\newcommand{\e}{\epsilon}
\newcommand{\Ga}{\Gamma}

\def\htwo #1 #2 {{\tilde{h}^{#1} _{#2}}}

\begin{document}
\title{An iterated random function with Lipschitz number one\footnote{This article was published in \emph{Theory of Probability and its Applications,} vol.~47 no.~2 (2003), pp.~286--300.}}

\author[Abrams]{Aaron Abrams}

\author[Landau]{Henry Landau}

\author[Landau]{Zeph Landau}

\author[Pommersheim]{James Pommersheim}

\author[Zaslow]{Eric Zaslow}

\date{\today}

\begin{abstract}
Consider the set of functions $f_{\th}(x)=|\th -x|$  on $\R$.  Define a 
Markov process that starts with a point $x_0 \in \R$ and continues with 
$x_{k+1}=f_{\th_{k+1}}(x_{k})$ with each $\th _{k+1}$ picked from a 
fixed bounded distribution $\mu$ on $\rp$.  We prove the conjecture of 
G. Letac that if $\mu$ is not supported on a lattice, then this process 
has a unique stationary distribution $\pi_{\mu}$ and any distribution 
converges under iteration to $\pi_{\mu}$ (in the weak-$^*$ topology).  
We also give a bound on the rate of convergence in the special case
that $\mu$ is supported on a two-point set.  We hope that the techniques 
will be useful for the study of other Markov processes where the transition 
functions have Lipschitz number one.
\end{abstract}

\maketitle

\section{Introduction}

In their recent paper \cite{persi}, P. Diaconis and 
D. Freedman describe a general 
method for studying Markov chains.  This consists
of viewing a state transition
of the chain as the action of a function chosen
at random from a family of 
state transformations; the Markov
process can then be expressed
as an iteration of these functions. 
When the state space is a metric space, e.g. $\R$,
they show that if the 
functions exhibit a certain average contractivity
(described in terms of Lipschitz numbers), then
running the iteration 
backward rather than forward produces geometric
convergence to a unique 
stationary distribution for the chain.
Many applications of this idea are developed.

An example of a Markov process not immediately
amenable to this method, 
originally proposed for study by Letac \cite{letac},
is defined on the non-negative reals
$\rp = [0, \infty)$ by 
starting with a point $x_0 \in \R^+$ and for $k=0,1,\ldots$,
setting 
\begin{equation} \label{irf}
x_{k+1}=f_{\th_{k+1}}(x_k)\equiv|\th_{k+1}-x_k|,
\end{equation} 
where the $\th_k \in  \rp $ are chosen independently
from a fixed distribution $\mu$ (on $\rp$).  
If we choose the starting point $x_0$ 
from a distribution $\pi _0$, then the above Markov process defines a 
distribution $\pi _n$ for $x_n$.

We will consider the case that $\mu$ has finite expectation;
it is known (e.g. \cite{feller}, p. 208) 
%
that in this case the distribution
\begin{equation} \label{stationary}
\pi_{\mu}(x)=\frac{1}{E(\th)}\int _0 ^x\Pr\{\th>y\}dy
\end{equation}
is a stationary probability distribution 
for this process (where $\th$ is chosen according
to the distribution $\mu$ and $E(\th)$ is its expectation).
G. Letac conjectured in \cite{letac} that if $\mu$
is not supported on
a lattice (that is, the set of integer multiples of a fixed
real number), then $\st $ is the unique
stationary distribution, and that the distributions $\pi_n$
of $x_n$ should converge to $\st$ in the weak-$^*$ topology.
Both these conjectures would follow from Theorem 1
of \cite{persi} if the $f_{\th}$'s had
Lipschitz numbers less than one (on average).


In this paper we prove Letac's conjecture for compactly supported $\mu$.

\begin{custom}{Theorem \ref{proof}}{ Let $\mu$ be a distribution on a 
bounded interval $Y\subset \rp$, with $\mu$ not supported on a lattice.
Then $\st$ 
is the unique stationary distribution for the process (\ref{irf}).  
Moreover, if $x_0$ is chosen 
according to an initial distribution $\pi_0$ then the distributions $\pi_n$
of $x_n$ converge to $\st$ in the weak-$^*$ topology.
}\end{custom}

We prove this theorem in Sections \ref{henrys} and \ref{backward}, and 
in Sections \ref{graphs} and \ref{boundsection} we 
analyze more carefully the special case of $\mu$ supported on just two
points.  In this case we provide a second independent proof of Letac's
conjecture, and we also obtain a bound on the rate of convergence of
the distributions $\pi_n$.  Our
hope is that the techniques will be useful for the study of 
other Markov processes where the transition functions have Lipschitz 
number one.  

After completing this work the authors learned of a paper by
J. Mattingly \cite{mattingly} addressing related questions; that paper
obtains results similar to ours.

Let $\mu$ be a probability measure on a space $Y$ and let 
$\{f_{\th}: \th \in Y\}$ be a set of maps from a metric space $X$ 
into itself.  The \emph{iterated random function induced by $\mu$} is
the Markov process of repeatedly picking 
$\th \in Y$ according to $\mu$ and then applying $f_{\th}$. 
For $x \in X$, and $\Th= ( \th_1, \th_2, 
\dots )$ with $\th_i \in Y$, let 
$B_n(\Th,x)$ be the $n$th \emph{backward iterate} 
of $x$ and $F_n(\Th,x)$ the $n$th \emph{forward iterate} of
$x$;  that is,
\begin{eqnarray*}
B_n(\Th,x)&=&f_{\th_1}\circ \cdots \circ f_{\th_n}(x), \\                                   
F_n(\Th,x)&=&f_{\th_n}\circ \cdots \circ f_{\th_1}(x).
\end{eqnarray*}
Letac observed that for fixed $n$ and variable $\Th$, the distributions of 
$B_n(\Th,x)$ and $F_n(\Th,x)$ are identical; yet for fixed $\Th$,
as $n$ tends to infinity the sequences $B_n(\Th,x)$ and $F_n(\Th,x)$
behave quite differently.  He used this to prove the following general
principle.

\begin{custom}{Letac's Principle \cite{letac}}{
If the backward iterates $B_n(\Th,x)$ converge for almost all $\Th$
to a limit $B(\Th)$ which is independent of the
starting point $x\in X$, then the distribution $\st$ on $X$
induced\footnote{For any set $S\subset X$, the induced probability 
distribution $\st(S)$ is the measure of the set $\{\Th\ |\ B(\Th)\in S\}$.}
by the map $\Th \mapsto B(\Th)$ is the unique
stationary distribution for the iterated random function
induced by $\mu$.
Moreover, if $x_0$ is chosen according to 
an initial distribution $\pi_0$, then the 
distributions $\pi_n$ of $x_n$
converge to $\st$ in the weak-$^*$ topology.
}\end{custom}

In Section \ref{henrys} we prove the following general result which
implies the hypotheses of Letac's principle.
Note that this theorem requires checking a property of 
only the \emph{forward} iterates.

\begin{custom}{Theorem \ref{reformulation}}{
Consider a probability measure $\mu$ on a space $Y$ and a set of functions 
$\{f_{\th} :\th\in Y \} $
on a complete metric space $X$.  Suppose 
there exists a subset $S\subset X$ such that 
\begin{itemize}
\item[(i)] $f_{\th}(S) \subset S$ almost surely for $\th $  picked from 
$\mu $, and
\item[(ii)] $\lim _{n \to \infty} \mbox{diam}(F_n(\Th , S)) =0$ almost 
surely, where $\Th=(\th_1,\th_2,\ldots)$ and the $\th_i$ are picked 
independently from $\mu$.
\end{itemize}
Then the backward iterates applied to an element of $S$ converge almost
surely and the limit is independent of the starting point (in $S$).
}
\end{custom}

In Section \ref{backward}
we establish the hypotheses of Theorem \ref{reformulation}, 
completing our proof of Letac's conjecture.

Our second proof of Letac's conjecture applies only in the special case
that $\mu$ is uniform on a two-point set $\{\a,\b\}$; however the
argument is independent of Letac's principle.
Our method is to analyze the iterated random function by studying the 
orbit of a point; we begin this study in Section \ref{backward} and
continue in more detail in Section \ref{graphs}.  This technique 
leads to a bound in Section \ref{boundsection} (Theorem \ref{zephthm}) 
on the rate of convergence of the distributions $\pi_n$ in this case.

\section{A General Result}\label{henrys}

In this section we prove a general theorem which lets us conclude that
the backward iterates converge by looking only at the forward iterates.

Let $\mu$ be a probability measure on a set $Y$ and $\{f_\th : \th \in Y \}$ 
a set of functions from a complete metric space $X$ into itself.  Denote by 
$\Th = (\th_1, \th_2, \dots )$ an element of 
$Y^{\infty} = Y \times Y \times \cdots$ 
and by $\mu ^{\infty}= \mu \times \mu \times \cdots$ the infinite product 
measure on $Y^{\infty}$.
The following lemma is used in the proof of Theorem \ref{reformulation} to 
move from considering backward iterates to considering forward iterates.

\begin{lem}\label{importantsimple}
Let $\Gamma \subset Y^{\infty}$, and let
$k_1, \dots, k_n$ be a finite set of distinct indices.  Let
$$\Gamma'= \{ \Th' : \mbox{there exists a } \Th \in \Gamma 
\mbox{ such that } \th'_{k_j}= \th _j , \ j= 1 \dots n \},$$
where $\Th=(\th_1, \th_2, \dots )$ and $\Th'= (\th'_1, \th'_2, \dots)$. 
Then $\mu^{\infty} (\Gamma) \leq \mu^{\infty} (\Gamma')$.
\end{lem}

\bpf  Let $\sigma$ be any finite permutation 
of $\N$ with $\sigma(j)=k_j$ 
for $j=1,\dots,n$.  (A finite permutation is one with $\sigma(j)=j$ for
all $j$ sufficiently large.)  Since $\mu ^{\infty}$ is invariant under the 
action of finite permutations, we have $\mu ^{\infty} (\Gamma)= \mu ^{\infty}
(\sigma(\Gamma))$.  But $\sigma(\Gamma) \subset \Gamma'$, so the result 
follows.
\epf

For $S\subset X$,
set $B_n(\Th , S)= \{ B_n(\Th , x) : x \in S \}$ and  
$F_n(\Th , S)=\{ F_n(\Th , x) : x \in S \}$.
Let $\mbox{diam}(S)= \sup _{x,y\in S} d(x,y) $ be the diameter of $S$
(where $d$ denotes the distance function on $X$).

\begin{thm} \label{reformulation}
Suppose there exists a subset $S$ of $X$ such that
\begin{itemize}
\item[(i)] $f_{\th}(S) \subset S$ for almost every $\th$ picked from $\mu$,
\item[(ii)] $\lim _{n \to \infty} \mbox{diam}(F_n(\Th , S)) =0$ for
almost every $\Th$ picked from $\mu^{\infty}$.
\end{itemize}
Then for all $x\in S$, the backward iterates $B_n(\Theta,x)$ converge 
for almost all $\Theta$, and the limit is independent of $x$ (in $S$).
\end{thm}

\bpf
Fix $\epsilon >0$.  Set $\Ga^i_N= \{ \Th: \mbox{diam}(F_n(\Th, S)) < 2^{-i} 
\mbox{ for all } n \geq N \}$.  Since $\Ga ^i_N \subset \Ga ^i _{N+1}$ and 
$\mu^{\infty}(\cup_n \Ga^i _n)=1$ (by hypothesis (ii)), there exists $N(i)$ 
such that $\mu^{\infty}(\Ga^i_{N(i)})> 1 - \epsilon 2^{-i}$.  Applying 
Lemma \ref{importantsimple} to the set $\Ga =\Ga^i_{N(i)}$ with 
$k_j= N(i)+1-j$, $j=1, \dots, N(i)$, we see that the measure of the set 
$\Lambda ^i_{N(i)}= \{ \Th: \mbox{diam}(B_{N(i)}(\Th, S))<2^{-i} \}$ is 
greater than or equal to $\mu^{\infty}(\Ga^i_{N(i)})>1 - \epsilon 2^{-i}$.  
Hypothesis (i) implies that for $n \geq N$, $B_n(\Th,S) \subset 
B_N(\Th, S)$ for almost all $\Th$; hence 
\begin{equation} \label{e:1} \mu^{\infty} \{\Th: \mbox{diam}(B_n(\Th,S))< 
2^{-i} \mbox{ and } B_n(\Th,S) \subset B_{N(i)}(\Th, S) \mbox{ for all } 
n \geq N(i) \} 
\end{equation}
\[= \mu^{\infty} ( \Lambda ^i_{N(i)})> 1 -\epsilon 2^{-i}. 
\]

Letting 
$$\Omega = \bigcap_i \{ \Th: \mbox{diam}(B_n(\Th,S))< 2^{-i} \mbox{ and } 
B_n(\Th,S) \subset B_{N(i)}(\Th, S) \mbox{ for all } n \geq N(i) \},$$
(\ref{e:1}) implies $\mu ^{\infty} (\Omega) \geq 1- \sum_i \epsilon 2^{-i}
=1-\epsilon$.  It follows from the definition of $\Omega$ that the 
$\{B_n(\Th, x) \}_n$ are Cauchy sequences (i.e. the backward iterates 
converge) for every $x \in S$, $\Th \in \Omega$.  Furthermore, since for 
$x,y \in S$, $|B_n(\Th,x) - B_n(\Th, y)| \leq \mbox{diam}(B_n(\Th,S))$, 
we see that for $\Th \in \Omega$, $\lim B_n(\Th,x)=\lim B_n(\Th,y)$.
Thus we have shown that on the set $\Omega$ of measure 
at least $1-\epsilon$, the backward iterates applied to an element of $S$ 
converge to a limit independent of the starting point.  Since $\epsilon$ 
is arbitrary, the result follows.\epf

\section{Proof of Letac's Conjecture}\label{backward}

We now return to the iterated random function described by equation
(\ref{irf}).
In this section, we consider $\mu$ supported on a bounded subset 
$Y=[0,b]$ of $\R ^+$.  Setting $S=[0,b']$ for any $b'\geq b$, 
it is clear that condition (i) of Theorem \ref{reformulation} is 
satisfied.  The rest of this section shows that $S$ and 
$\{f_{\th} : \th \in [0,b] \}$ satisfy condition (ii) of Theorem 
\ref{reformulation}, and then applies the conclusion of Theorem 
\ref{reformulation} to complete the proof of Theorem \ref{proof}.  

The argument showing that $S$ and $\{f_{\th} : \th \in [0,b] \}$ 
satisfy condition (ii) of Theorem \ref{reformulation} can be summarized 
as follows:  on an interval $I$, $f_\th$ acts in one of two ways.
If $\th \not\in I$, then $f_{\th}$ translates $I$; whereas if 
$\th \in I$, then $f_{\th}$ translates and ``folds'' (and in particular
shrinks) $I$.  Furthermore, for $\th$ near the center of $I$, $f_\th$ 
will shrink $I$ by about half.  Condition (ii) follows 
from showing that the process defined by equation $(1)$ (for suitable 
$\mu$) almost halves an interval infinitely many times.
The key fact we will use is that the orbit of a point is dense.

\subsection{The orbit of a point}
 
In this section we discuss the orbit of a point $x$ under
iterates of the two functions $f_{\a}(x)=|\a-x|$ and $f_{\b}(x)=|\b-x|$,
where $0<\a<\b$.
We use the standard notation $x \bmod y$ to mean the unique number $z$,
$0\leq z < y$ such that $z-x$ is an integer multiple of $y$.  
Let $x\in[0,\b]$.

The {\em orbit
${\mathcal O}_x$ of $x$} is the set of points $y$
such that
$y=f_{\th_n}\circ\cdots\circ f_{\th_1}(x)$ for some $n,$
with each $\th_i \in\{\a,\b \}$.
 
\begin{lem}\label{orbits}  For each $x\in[0,\b]$, we have ${\mathcal O}_x
=\left\{ n\a+\e x \bmod \b : n\in\Z, \e=\pm 1\right\}$.
\end{lem}
 
\bpf
Let $S_x=\left\{ n\a+\e x \bmod \b : n\in\Z,
\e=\pm 1\right\}$.
As shorthand, we represent the real number $ n\a+\e x \bmod \b$
by the ordered pair $(n,\e)\in\Z\times\{\pm1\}$.
Note that $f_{\b}((n,\e))=\b- (n\a+\e x \bmod \b) =(-n,-\e)$.
Also, $$f_{\a}((n,\e))=
\begin{cases}
(n-1,\e) & \text{if $(n,\e)\geq\a$,} \\
\a-(n\a+\e x\bmod\b)=(-(n-1), -\e) & \text{if $(n,\e)<\a$;}
\end{cases}$$
thus the set $S_x$ is closed under $f_{\th}$ for each                           
$\th\in{\{\a,\b \}}$
and we have ${\mathcal O}_x\subseteq S_x$.
 
To show $S_x\subseteq {\mathcal O}_x$, it suffices
to show that
\begin{enumerate}
\item{} $(n,\e)\in {\mathcal O}_x
\implies (-n,-\e)\in{\mathcal O}_x$, and
\item{} $(n,\e)\in {\mathcal O}_x
\implies (n-1,\e)\in{\mathcal O}_x$.
\end{enumerate}
 
The first claim follows by applying $f_{\b}$.
The second follows by
applying $f_{\a}$ if $n\a+\e x\bmod\b\geq \a$,
or by first applying
$f_{\a}$ and then applying $f_{\b}$ if $n\a+\e x\bmod\b<\a$.
\epf
 
\begin{cor} \label{c1}  Let $0<\a<\b$.  If $\a/\b$ is irrational
then the orbit of any point under iteration of the functions 
$f_{\a}, f_{\b}$
is dense in $[0,\b]$.  If $\a/\b=p/q$ with $p,q\in\N$ relatively
prime, then the orbit of any point intersects each closed subinterval
of $[0,\b]$ of length at least $\b/q$.
\end{cor}

\bpf  If $\a/\b$ is irrational then the set $\{n\a\bmod\b\}$ is
dense in $[0,\b]$, so the orbit of a point is also dense in $[0,\b]$.  
If $\a/\b=p/q$ 
then the orbit of $x$ contains the set $\{n\a+\e x\bmod\b\}=
\{\frac kq \b+x \bmod\b |k=0,\dots,q-1\}$ which intersects 
each closed subinterval of $[0,\b]$ of length at least $\b/q$.
\epf

\subsection{Shrinking intervals}

Denote the length of an interval $I$ by $\length{I}$.

The following lemma gives a criterion for recognizing when the interval 
$f_{\theta_n} \circ \cdots \circ f_{\theta_1}(I)$ is no bigger than about 
half the size of $I$.

\begin{lem} \label{shrinkage}
Let $\Th=(\th_1, \th_2 \dots)$, with $\th_i \in \R^+$, and let $I\subset \rp$ 
be a bounded interval with midpoint $m$.  Fix $\d_1, \d_2 > 0$.
If some $x$ with $|x-m|<\d_1$ satisfies $F_n(\Th,x)<\d_2$, then 
$$\length{F_n(\Th, I)}\leq \frac{1}{2}\length{I}+\d_1+\d_2.$$
\end{lem}

\bpf
The functions $f_{\th}$ have Lipschitz number 1, so for all 
$y\in I$, we have
$$|F_n(\Th,x)-F_n(\Th,y)|\leq|x-y|\leq \frac{1}{2}\length{I}+\d_1.$$
Since $F_n(\Th,x)\leq\d_2$, it follows that $F_n(\Th,y)\leq
\frac 1 2 \length I + \d_1+\d_2$.  Also $F_n(\Th,y)\geq 0$, so
the result follows.\epf

We say that $x\in\rp$ is a {\em point of increase}
for the probability distribution $\mu$ if $\mu(I)>0$
for every open interval $I$ containing $x$.  The nonzero
real numbers
$\a$ and $\b$ are {\em irrationally related} if
$\a/\b$ is irrational.

\begin{lem} \label{backbone}
Let $\mu$ be a probability distribution supported on
a bounded interval
$Y \subset \rp$.  Suppose that $\mu$ is not supported on a lattice.
%
%
%
%
Then for any interval $I\subset \R^+$ and any $\epsilon >0$, there
exist $N=N(I,\e)$ and intervals $J_1, \dots, J_N$
with $\mu (J_i) > 0$, such that for all $(\theta_1,
\dots, \theta_N) \in J_1 \times
\dots \times J_N$,
\[ \length {f_{\th_N}\circ\dots\circ f_{\th_1}(I)} < \epsilon \]

\end{lem}

\bpf
We first assert that for any interval $K$, and
for any $r>0$,
there exist intervals $J_1, \dots J_n$ with
$\mu (J_i) >0$ such that for all $(\theta_1,
\dots, \theta_n) \in J_1 
\times \dots \times J_n$,
\[ \length {f_{\th_n}\circ\dots\circ f_{\th_1}(I)} < 
(\frac{1}{2}+r)\length {K} . \]
The result then follows by repeating this shrinking
process.

To prove the assertion,
we will first construct the midpoints $m_i$ of the intervals $J_i$.
We consider two cases;
as $\mu$ is not supported on a lattice, either (a) $\mu$ has two irrationally
related points of increase, or (b) $\mu$ has an infinite set of points 
of increase which are pairwise rationally related. 

Let $m$ be the midpoint of the interval $K$.

First, suppose (a) holds.
Let $\alpha<\beta$ be two irrationally related points of increase
for $\mu$.  By Corollary \ref{c1} the orbit of $m$ 
under iteration of the functions $f_{\a}$ and $f_{\b}$ is
dense in $[0,\b]$; thus 
there exists a sequence
$m_1, \dots, m_n$ with each
$m_i$ equal to $\a$ or $\b$ such that
\begin{equation}\label{shrunk}
f_{m_n} \circ \dots \circ f_{m_1}(m) < \frac{r}{2} \length {K} .
\end{equation}

Now suppose (b) holds.
In this case we can find points of increase $\a,\b$ of $\mu$ 
with $\a/\b=p/q$, with $p,q\in\N$ relatively prime and $q$ arbitrarily
large.  (Note that we cannot simultaneously guarantee that $\a<\b$.)
Find such $\a,\b$ with $q>2\b/(r\length{K})$.  By Corollary \ref{c1},
the orbit of $m$ under the functions $f_{\a}$ and $f_{\b}$
intersects the closed interval $[0,\b/q]\subset[0,r\length{K}/2]$.
Thus again there exist $m_1,\dots, m_n$ with each
$m_i$ equal to $\alpha$ or $\beta$ such that (\ref{shrunk}) holds.

We have now defined the $m_i$ in each case.  Set $J_i$ to be an interval 
with midpoint $m_i$ and length no bigger than
$\frac{r}{2n} \length{K}$.  Note that $J_i$ has positive
measure because every $m_i$ is a point of
increase of $\mu$.

To verify the claim, note that for any
$(\theta_1, \dots, \theta_n)
\in J_1 \times \dots \times J_n$, we have
\[f_{\th_n} \circ \dots \circ f_{\th_1}(m)
< \frac{r}{2}\length {K} +  \frac{nr}{2n}\length {K}
=r\length {K} ,\] 
since each $f_{\th}$ has Lipschitz number 1.
Thus we can apply the previous lemma with $x=m$, any $\d_1>0$,
and $\d_2=r\length K$.

\epf

\begin{thm} \label{proof}
Let $\mu$ be a distribution on a bounded interval $Y\subset \rp$, with
$\mu$ not supported on a lattice.
Then $\st$ is the unique stationary 
distribution for the process (\ref{irf}).  
Moreover, if $x_0$ is chosen according to an initial 
distribution $\pi_0$ then the distributions $\pi_n$
of $x_n$ converge to $\st$ in the weak-$^*$ topology.
\end{thm}
\bpf
We will show that $\{f_{\th} : \th \in Y\}$, $S=[0,b']$, $b'\geq b$ 
satisfy conditions (i) and (ii) of Theorem \ref{reformulation}.  
Consequently,  the backward iterates converge independent of a 
starting point chosen in $\bigcup _{b' \geq b}[0,b']= \R ^+$ 
and the application of Letac's principle completes the proof.

It is clear that $S$ 
satisfies condition (i) of Theorem \ref{reformulation} for all 
$\Th \in Y^{\infty}$.  To establish (ii), fix $\epsilon>0$.  By
Lemma \ref{backbone} (with $I=[0,b']$), there exist
intervals $J_1, \dots, J_N$ such that 
$$ \length { f_{\th_N}\circ \cdots
\circ f_{\theta_1}(I) } < \epsilon$$
for all 
$(\th_1, \dots \th_N) \in J_1 \times \dots \times J_N$.
As the $f_{\th}$'s do not lengthen intervals, we have 
$\length{ F(\Th,S) } < \epsilon$ 
for any $\Th=(\th_1,\th_2,\ldots)$ such that for some 
$K$, $\th_{K+i}\in J_i$ for each $i=1,\ldots,N$.
Since each $\mu(J_i)>0$, this occurs with probability 1.  Since $\epsilon$
was arbitrary, condition (ii) of Theorem 
\ref{reformulation} is satisfied. \epf

\section{Special Case:  Two-point Distributions} \label{graphs}

Our second proof of Theorem \ref{proof} relies on a detailed
understanding of the orbits $\O_x$.  In this section 
suppose $\th$ is chosen from the uniform distribution
$\mu$ on the set $\{\a,\beta\}$ with $0<\a<\beta;$
i.e. ${\rm Pr}(\th= \a) = {\rm Pr}(\th=\b) = 1/2$.  By scaling, 
we may assume $\beta=1$; we do so for the remainder of the paper.
We use the notation $\langle x\rangle$ to denote the number $x\bmod 1$. 

\subsection{The graph of an orbit}
\label{graphoforbit}

Recall Lemma \ref{orbits}, and note that 
the proof of Lemma \ref{orbits} (with $0<\a<\b=1$) reveals the complete
structure of the orbit of $x \leq 1$.  From now on assume $\a$ is
irrational.
Notice that the numbers $\langle n\a+\e x\rangle$ are
distinct unless $x$ is in
the orbit of $0$ (which is also the orbit
of $1$ and $\a$), the orbit 
of $1/2$, the orbit of $\a/2$, or the
orbit of $\frac{1+\a}{2}$.  These
four orbits we call \emph{singular}; other
orbits are \emph{generic}.

To visualize the orbit $\O_x$, we associate to each point $x$ a labeled
directed graph $G_x$, as follows.  The vertex set of $G_x$ is $\O_x$, 
with the vertex $\langle n\a+\e x\rangle$ labeled with the ordered pair
$(n,\e)$.  This label is unique if $\O_x$ is generic.  A directed edge 
goes from $(n,\e)$ to $(n',\e')$ if there exists a $\th$ such that 
$f_{\th}(\langle n\a+\e x\rangle)=\langle n'\a+\e' x\rangle$; the proof
of Lemma \ref{orbits} tells us exactly where to put these edges.  

Note that
if $x'\in\O_x$ then $G_{x'}$ looks exactly like $G_x$ but with different 
labels.  Also note that if $\O_x$ is singular then a vertex of $G_x$ will
have multiple labels due to coincidences among the numbers 
$\langle n\a+\e x\rangle$; e.g. a vertex labeled $(n,\e)$ in $G_0$ is
also labeled $(n,-\e)$.

We illustrate the case $\a=1/\sqrt{2}$.  Figure 1(a) shows the (generic)
graph $G_{.2}$, and Figure 1(b) shows the (singular) graph $G_0$.  In
the singular orbit we have consolidated the pair of labels $(n,\e)$
and $(n,-\e)$ into the single label $n$.  Both graphs have been drawn
so that the vertical edges represent an application of $f_1$ and the
horizontal and diagonal edges represent an application of $f_{\a}$.
An undirected edge indicates two oppositely oriented edges.  
Note that each vertex has out-degree two, since $\th$ has two 
possible values.

\begin{figure}[h]
\begin{center}
\scalebox{.5}{\includegraphics{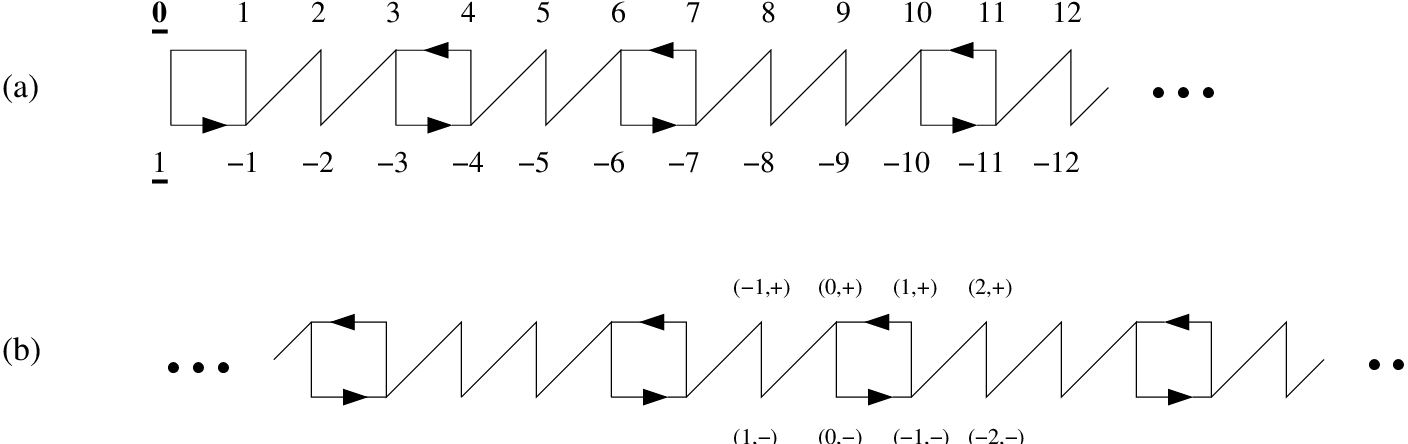}}
\caption{Two orbits for $\a=1/\sqrt{2}>1/2$.  
The orbit of $x=0.2$ is generic (a), and the orbit of $x=0$ is singular (b).  
In the orbit of $0$, vertices are labeled only with $n$, since 
$\e$ doesn't matter.  The boldface labels
\underbar{{\bf 0}} and \underbar{{\bf 1}} refer to
the actual numbers $0$ and $1$.}
\label{largealpha}
\end{center}
\end{figure}

The large-scale properties of the graphs $G_x$ (for general irrational
$\a$) are important as well.  Suppose for a moment that $\a>1/2$.
There are three types of vertices in a generic graph (see Figure 1):
\begin{itemize}
\item{} \emph{small}, where $\langle n\a + \e x \rangle < 1-\a$;
\item{} \emph{medium}, where $1-\a < \langle n\a + \e x \rangle < \a$; and
\item{} \emph{large}, where $\langle n\a + \e x \rangle > \a$.
\end{itemize}
The large vertices have in-degree one, and appear in $G_x$ as upper-right
and lower-left corners of boxes.  The small vertices have in-degree
three, and form the other corners of boxes.  The medium vertices
are the remaining ones, which are not part of any box.  The fact
that $\langle n\a\rangle$ is uniformly distributed
mod 1 implies that
the ratio of the numbers of small : medium : large vertices is
$1-\a:2 \a-1:1-\a$.
If we write $\a=q(1-\a)+r$ with $q$ a positive integer and $0<r<1-\a$, then
the number of diagonal edges between consecutive boxes is either 
$q$ or $q+1$. 
These numbers occur in a pattern depending on $x$ but always in
the ratio $(1-\a-r) : r$.  For $\a=
1/\sqrt{2}$ (Figure \ref{largealpha}) we see that there are either 
2 or 3 diagonal edges between consecutive boxes, this number being 2 
about $58.6\%$ of the time (since 2 occurs exactly
$\frac{3-4\a}{3\a-2}=\sqrt{2}$ times as often 3).  

\begin{figure}[h]
\begin{center}
\scalebox{.5}{\includegraphics{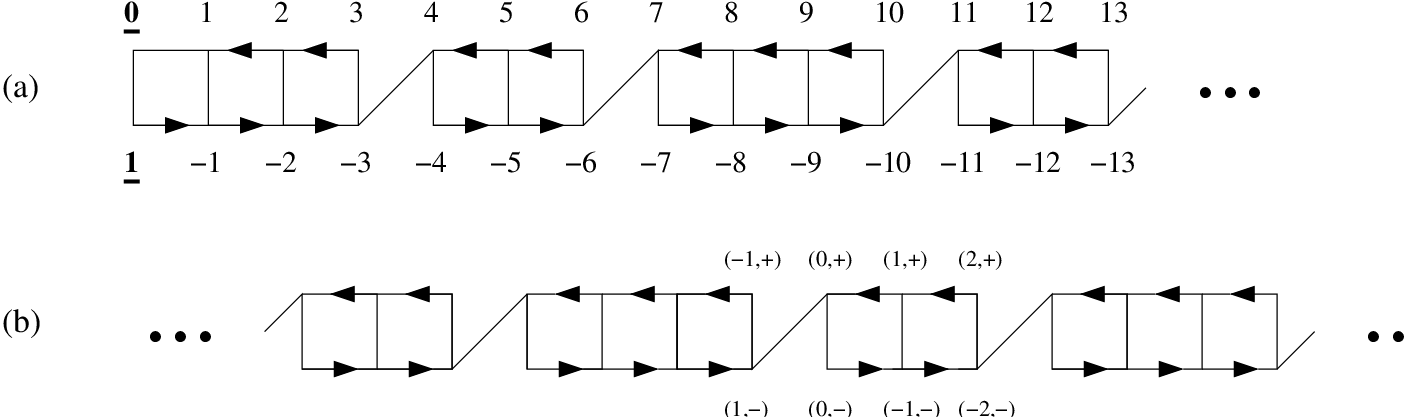}}
\caption{Two orbits for $\a=1-1/\sqrt{2}<1/2$.  Again,
the orbit of $x=0.2$ is generic (a) and the orbit of  $x=0$ is singular (b).}
\label{smallalpha}
\end{center}
\end{figure}

In the case $\a<1/2$, the graphs are similar, but we exchange the roles
of $\a$ and $1-\a$ in the definitions of small, medium, and large.
The graphs $G_x$ in this case have more boxes and fewer diagonal edges.
See Figure \ref{smallalpha}, where $\a=1-1/\sqrt{2}$.  The large vertices 
are the 
upper-right and lower-left 
corners of the rectangular strings of boxes; the small vertices are the
other corners of the rectangles; and the medium vertices are those 
vertices which are part of two distinct boxes.  The ratio of numbers of
small : medium : large vertices is $\a:1-2\a:\a$.
This time writing $1-\a=q\a+r$,
the number of boxes between consecutive diagonal edges is again
either
$q$ or $q+1$,
and these occur in the ratio $(\a-r) : r$.

Incidentally if $\a$ is rational then the graphs of orbits are finite,
but they exhibit many properties analogous to the infinite case.


\subsection{Uniqueness of stationary distribution}

Assume $\a\in(0,1)$ is irrational and $\mu$ is uniform on $\{\a,1\}$.
By (\ref{stationary}), the stationary
distribution $\pi_{\mu}$ is given by
\begin{eqnarray}\label{stationaryformula}
\st (x)= \begin{cases}
2x/(1+\a) &\mbox{ if  } 0\leq x\leq\a \\
(x+\a)/(1+\a) &\mbox{ if  } \a\leq x\leq 1.
\end{cases}
\end{eqnarray}

Let $\pi$ be any stationary distribution for the iterated random 
function (\ref{irf}) induced by $\mu$.  We will think
of $\pi$ as its distribution function, namely a right-continuous
increasing function $\pi:[0,1]\to[0,1]$; thus $\pi(x)$ is the amount
of mass which is concentrated on $\{y : y\leq x\}$.

We will show that the values of $\pi$ are determined by $\a$; the 
values we compute agree with 
the corresponding values for $\st$.

Being right-continuous, $\pi$ is determined by
its values on a dense set; with the aid of Figures \ref{largealpha}(a)
and \ref{smallalpha}(a)
we shall calculate $\pi(x)$ for $x\in\O_0=\{\langle n\a\rangle\}$.

By stationarity, the mass of $\pi$ is concentrated on $[0,1]$;
thus $\pi(1)=1$.  
We begin by establishing that $\pi$ cannot have
an atom (that is, a jump discontinuity) at $\langle n \a \rangle$ 
for any $n \in  \Z$.

\begin{lem}
\label{noatom}
Let $0<\a<1$ be irrational, let $\mu$ be uniform on $\{\a,1\}$,
and let $\pi$ be a stationary distribution for the associated 
iterated random function.  Then
$\pi(0)=0$ and $\pi$ is continuous at
$x=\langle n \a \rangle$, $n= \pm 1, \pm 2, \dots$.
\end{lem}

\bpf 
Suppose $1/2<\a<1$.  If a distribution $\pi_0$ has an atom of 
mass $m$ at $0$, then (referring to the graph $G_0$) we see that 
$\pi_1$ will have
an atom of mass $m/2$ at $1$, since the only way for $x_1$ 
to equal
$1$ is if $x_0=0$ and $\th_1=1$.  Thus if $\pi=\pi_0$ is stationary
and $\pi$ has an atom of mass $m$ at $0$, then $\pi$ also has an atom 
of mass $m/2$ at $1$.

Likewise, iterating the random function shows that at each vertex of 
$G_0$ there must
be an atom having mass equal to half the sum
of the masses
of the atoms at all incoming vertices.  It follows that
there must
be an atom of mass $3m/2$ at $\a$,
an atom of mass $2m$ at $\langle-\a\rangle$, and thereafter 
an atom of mass $2m$ at every small and medium vertex and 
an atom of mass $m$ at every large vertex.

As the
total mass of a distribution must be 1, it is impossible to
have $m>0$.

Thus $\pi(0)=0$.  We now note that if
$\pi$ were not continuous at
some point $\langle n\a \rangle$, then $\pi$ would have an
atom there;
then by the same reasoning $\pi$ would have an atom at $0$,
a contradiction.

The case $0<\a<1/2$ is similar (see e.g. Figure \ref{smallalpha}(a)).
Here an atom of mass $m$ at $0$ implies an atom 
of mass $m/2$ at $1$,
an atom of mass $3m/2$ at $\a$, and an atom 
of mass $m$ at $1-\a$,
and thereafter atoms of mass $2m$ at each small vertex 
and of mass
$m$ at each medium and large vertex.  As before, $m$ must be $0$, and
the result follows.
\epf

\medskip

We are now able to compute the values of $\pi(x)$ for
$x\in{\mathcal O}_0$,
using a technique similar to that employed in the lemma. 
Observe that the stationarity of $\pi$ is equivalent to the condition 
that for each $x\in[0,1]$ we have 
\begin{eqnarray}\label{dist}
\pi(x)&=&\sum\limits_{\th\in{\{1, \alpha \}}}\mu(\th)
\left(\sum\limits_{ \{y:f_{\th}(y)=x\}}\pi (y)\right) \\
&=&\frac{1}{2}\left(\pi(1+x)-\pi_-(1-x)\right)
+ \frac{1}{2} \left(\pi(\a+x)-\pi_-(\a-x)\right), \nonumber
\end{eqnarray}
where $\pi_-(y) \equiv \lim_{z\rightarrow y^-}\pi(z)$
(the limit exists since $\pi$ is increasing).  We will only need to 
consider the points $\langle n\alpha\rangle$; Lemma \ref{noatom} 
establishes that $\pi$ and $\pi_-$ are equal at these points.
Recall, of course, that $\pi(x)=1$ for all
$x\geq 1$ and $\pi(x)=0$
for all $x< 0$.


Again we first treat the case $1/2<\a<1$.
We already know $\pi(0)=0$ by Lemma \ref{noatom},
and also $\pi(1)=1.$ 
Suppose $$\pi(\a)=z\in[0,1].$$  Evaluating (\ref{dist}) at $x=\a$,
we get $\pi(\langle -\a\rangle)=2-2z$.  Similarly, at $x=\langle -\a\rangle$ we 
get $\pi(\langle 2\a\rangle)=3z-2$,
and so on. 
Using (\ref{dist}) repeatedly, we can 
compute $\pi(\langle n\a\rangle)$
for all $n\in\Z$.
In particular, if $n>0$ and $\langle n\a\rangle$ is small,
then we get $\pi(\langle n\a\rangle)=
(2n-L_n)z-(2n-2L_n)$ where $L_n$ is the number
of large vertices among $\langle \a \rangle, \ldots, \langle
n\a \rangle.$
Recall that $L_n/n \to 1-\a.$ 

We now claim that there is only one possible value of $z$. 
To see this, choose an increasing sequence
of positive integers $\{ n_i\}$ such that
$\langle n_i\a\rangle\to 0$.
We may assume the $\langle n_i\a\rangle$ are small.
Since $\pi(0)=0$ and $\pi$ is right
continuous, we have
$0=\lim_{i\to \infty}\pi(\langle n_i\a\rangle)=\lim_{i\to\infty}
[(2n_i-L_{n_i})z-(2n_i-2L_{n_i})]$,
whence $$z=\lim_{i\to\infty}\frac{2n_i-2L_{n_i}}{2n_i-L_{n_i}} 
= \frac{2\a}{1+\a}.$$

If $0<\a<1/2$ we proceed similarly and find that for 
$n>0$ and small 
$\langle n\a\rangle$, $\pi(\langle n\a\rangle)=
(2n-LM_n)z-(2n-2LM_n)$ where $LM_n$ is the number of 
large and medium vertices among $\langle \a \rangle, 
\ldots, \langle n\a \rangle.$  
We have $LM_n/n \to 1-\a$, and the corresponding limiting
argument again shows that $z=2\a/(1+\a)$.

Note that by (\ref{stationaryformula}), $\pi_{\mu}(\a)=2\a/(1+\a)$.

\begin{thm} \label{t1}
Let $\mu$ be uniform on the set $\{\a,1\}$, with
$0<\a<1$ and $\a$ irrational.  Let $\pi$ be a stationary
distribution for the iterated random function.
Then $\pi=\st$.
\end{thm}

\bpf
Being right continuous, $\pi$ is
determined by its values
on the dense set ${\mathcal O}_0$. 
Therefore, by the above discussion,
$\pi$ is determined by its value at $\a$. 
Finally, we have shown that
there is only one possible value of $\pi(\a)$,
namely the value $\st (\a)$.
\epf

\section{Rate of convergence}
\label{boundsection}

We will use  the description of $G_x$ to give us a rate of convergence 
for the backwards iterates in Theorem \ref{zephthm}.  We shall need 
some preliminary results before proving the bound in Section \ref{bound}.

\subsection{Continued fractions}
\label{contfrac}
Fix $\a$ irrational, and let
\begin{equation}\label{definea_i}
\a = a_0 +\frac{1}{\displaystyle a_1 +
\frac{1}{\displaystyle a_2 + \frac{1}{\ddots}}} 
\end{equation}
be the continued fraction expansion of $\a$ (each
$a_i$ is an integer, with $a_i>0$ if $i>0$).
For $n\geq 0$ let $p_n$ and $q_n$ be relatively prime
integers such that 
\[ \frac{p_n}{q_n} = a_0 +\frac{1}{\displaystyle a_1
+ \frac{1}{\displaystyle a_2 + \frac{1}{\ddots +
\frac{1}{\displaystyle a_n}}}}. \]
It follows that $q_n$ satisfies the recursion
$q_{n+1}= a_n q_n + q_{n-1}$ (see e.g. \cite{hardywright}).  
The $p_n$ and $q_n$ are such that 
\begin{equation} \label{e:approx}
|\alpha - \frac{p_n}{q_n}| < \frac{1}{2q_n^2}.
\end{equation}

We shall make use of the following fact.

\begin{lem} \label{lemmaclose}  For any $x\in[0,1]$, there is an
integer $k$ with $0\leq k<q_n$ such that 
$0\leq\langle x-k\a \rangle <\frac{3}{2q_n}$. 
\end{lem}

\bpf
Let $S=\{\langle k\a\rangle\ :\ k=0,1,\ldots,q_n-1\}$, and for each $i=
0,1,\ldots,q_n-1$ let $\ibar$ denote the unique integer among
$0,1,\ldots,q_n-1$ satisfying $\ibar p_n\equiv i (\bmod q_n)$.
Note that $0\in S$.  It essentially follows from (\ref{e:approx})
that each of the $q_n-1$ (disjoint) intervals
$(\frac i{q_n}-\frac 1{2q_n},\frac i{q_n}+\frac 1{2q_n})$, where 
$i=1,\ldots,q_n-1$, contains exactly one point of $S$, namely 
$\langle \ibar a\rangle$.  To see this, note that by (\ref{e:approx}),
$$\left|\ibar \a-\ibar \frac{p_n}{q_n}\right|=\ibar\left|\a-
\frac{p_n}{q_n}\right|<\ibar\cdot\frac 1{2q_n^2}<\frac{1}{2q_n}.$$
This implies that there is no integer strictly between $\ibar\a$ and
$\ibar p_n/q_n$, and therefore $\left|\ibar\a-\ibar p_n/q_n\right|=
\left|\langle\ibar \a\rangle-\langle\ibar p_n/q_n\rangle\right|$.  Thus
$$\left|\langle\ibar \a\rangle-i/q_n\right|
=\left|\langle\ibar \a\rangle-\langle\ibar p_n/q_n\rangle\right|
=\left|\ibar\a-\ibar p_n/q_n\right|<\frac1{2q_n}.$$

Note also that if $\a>p_n/q_n$ then $\langle\ibar \a\rangle>i/q_n$ for 
all $i=0,1,\ldots,q_n-1$; whereas if $\a<p_n/q_n$ then $\langle\ibar 
\a\rangle<i/q_n$ for all such $i$.  Thus no two consecutive points of 
$S$ are more than $\frac{3}{2q_n}$ apart.

To summarize, the smallest point of $S$ is $0$, any two consecutive 
points of $S$ are within $\frac{3}{2q_n}$ of each other, and the largest
point of $S$ is greater than $\frac{q_n-1}{q_n}-\frac{1}{2q_n}=
1-\frac{3}{2q_n}$.  This proves the lemma.
\epf

\subsection{Random walks on the line}
\label{rwonline}

Fix $0<\alpha<1$ irrational, and fix $0<x_0<\alpha$ in a generic orbit. 
The vertex $v_0$ corresponding to $x_0$ in the graph $G=G_{x_0}$ 
is small; thus (see Section \ref{graphoforbit}) $v_0$ is a cut-vertex of $G$
(that is, $G-v_0$ is not connected).  

For an arbitrary vertex $v$ of $G$, let $d(v_0,v)$ be the 
distance from $v_0$ to $v$ in $G$, ignoring orientations.
Define a map $\rho$ from the vertex set of $G$ to $\Z$ 
as follows:  let $\rho(v)=n>0$ if $d(v_0,v)=n$ and $v$ is in the 
right hand component of $G-v_0$; let $\rho(v)=n<0$ if $d(v_0,v)=-n$
and $v$ is in the left hand component of $G-v_0$; and let
$\rho(v_0)=0$.

\begin{custom}{Observation}
If $\th$ is chosen uniformly from $\{\a,1\}$ then for any $x\in\O_{x_0}$,
we have $$\rho(f_{\th}(x))=
\left\{ \begin{array}{r@{\quad}l}
\rho(x)+1 & \mbox{w.p. }1/2,\\
\rho(x)-1 & \mbox{w.p. }1/2.
\end{array}\right.$$
\end{custom}

The next lemma says that very long paths on $G$ must go through vertices
which correspond to very small numbers.  As $x_0$ is fixed, we no 
longer distinguish between a point $x\in\O_{x_0}$ and the corresponding
vertex of $G$.

\begin{lem} \label{l:farsmall}
Let $x,y \in \O_{x_0}$, and suppose that $\rho(x)-\rho(y)\geq 2q_m$.  Then 
there exists $z\in\O_{x_0}$ with $\rho(x)\geq\rho(z)\geq\rho(y)$ such
that $0\leq z<\frac{3}{2q_m}$.
\end{lem}

\bpf
From the graph $G$ (more precisely, from the proof of Lemma \ref{orbits}),
it is clear that for any $k\in\Z$,
$$k\leq|\rho(\langle x-k\a\rangle)-\rho(x)|\leq 2k.$$  

Suppose first that $x=\langle r\a+x_0\rangle$ for some $r\in\Z$;
thus the vertex $x$ is on the top row of the graph $G$ and $\langle
x+\a\rangle$ is to the right of $x$.
Lemma \ref{lemmaclose} gives a number $k$ with $0\leq k < q_m$ and
$0\leq\langle x-k\a\rangle:=z<\frac{3}{2q_m}$. 
Then $\rho(y)\leq\rho(x)-2q_m\leq \rho(z)\leq\rho(x)$
as desired.  

On the other hand if $x$ is on the bottom row of $G$, i.e. $x=\langle
r\a-x_0\rangle$ for some $r\in\Z$, then we apply Lemma \ref{lemmaclose}
to the number $\langle x+(q_m-1)\a\rangle$ to get a $k$ with $0\leq k<q_m$
and $0\leq\langle x+k\a\rangle:=z<\frac{3}{2q_m}$.  Again we have
$\rho(y)\leq\rho(x)-2q_m\leq\rho(z)\leq\rho(x)$.
\epf

We shall need the following; the simple proof provided was pointed out 
to us by D. Aldous.

\begin{lem} \label{l:randwalk}
There exists an integer $N_0$ and a constant $c_0>0$ such that for 
all $n>N_0$, a simple 
random walk of length $n^3$ on a line
remains within $n$ units of its starting point
with probability less than or equal to
 $e^{-c_0n}$.
\end{lem}
\bpf   For a walk to stay within $n$ units of its starting point
it is necessary (though not sufficient) that no subwalk of the walk moves
more than $2n$ units in either direction.
By the central limit theorem, there is an integer $N_0'$ and a
constant $c_0'>0$ so that for $k>N_0'$, a walk of length $c_0'k^2$ 
remains within $k$ units of its starting point with probability 
less than $e^{-1}$.
Given a walk of length $n^3$, we break it into pieces of length $c_0'(4n^2)$
and apply the previous fact with $k=2n$ to each of the $n/4c_0'$ pieces.
It follows that the probability that the walk
remains within $n$ units of its starting point is at most
$e^{-n/4c_0'}$.
\epf

\subsection{The bound}
\label{bound}

The content of the following theorem is essentially that with
high probability, the entire interval $[0,1]$ gets within $\e$ 
of its limit under backward iteration after $O((1/\e)^3
\log(1/\e))$ steps.

\begin{thm}
\label{zephthm}
Let $0<\a<1$ be irrational, and fix $\e$ small and positive.
Find $k$ such that 
$8/q_k<\e$ (where $q_k$ is as in Section \ref{contfrac}), and let 
$N=8q_k^3\log_2q_k$.  There is a constant $c>0$, which does not depend 
on $\e$, such that with probability larger than $1-e^{-c/\e}$, we have
\begin{equation}\label{conclusion}
|B_n(\Th,x)-B(\Th) | < \e \quad\mbox{for all $n\geq N$ and for all 
$x\in[0,1]$,}
\end{equation} 
where $\Th=(\th_1,\th_2,\ldots)$ is chosen uniformly from 
$\{\a,1\}^{\infty}$.
\end{thm}

\bpf
Let $j=8q_k^3$.  Given an interval $I$ we compute a lower
bound $\mathsf{p}$ for the probability that
\begin{equation}\label{intervalshrink}
\length{f_{\theta_1}\circ \cdots \circ f_{\theta_j}(I)}
\leq \frac{1}{2}\length{I} + \frac{2}{q_k}. 
\end{equation}
It follows then that with probability bigger than $\mathsf{p}^t$,
\[\length{f_{\theta_1}\circ \cdots \circ
f_{\theta_{tj}}([0,1])} \leq (\frac{1}{2})^t 
+ \frac{4}{q_k}. \]  
In this way, with 
$t= \lfloor\log_{\frac{1}{2}} \frac{4}{q_k} +1\rfloor= 
\lfloor\log_2 q_k -1 \rfloor $
and $n \geq tj$ fixed, we have that for all $x\in [0,1]$, 
(\ref{conclusion}) holds, i.e.,
\begin{eqnarray*}
\left|B_n(\Th,x)-B(\Th)\right|&=&\lim_{m \to \infty} 
|B_n(\Th,x) - f_{\theta_1}\circ \cdots \circ f_{\theta_m}(x)| \\
&=& \lim_{m\to\infty}|B_n(\Th,x)-B_n(\Th,f_{\th_{n+1}}\circ\cdots\circ
f_{\th_m}(x)| \\
&<& \frac{8}{q_k} <\e,\\
\end{eqnarray*}
with probability at least $\mathsf{p}^t$.

We get the estimate $\mathsf{p}$  as follows.  
Pick a point $m^* \in \O_{x_0}$ ($x_0$ as in Section 
\ref{rwonline})
such that $m^*$ is within $\delta=\frac{1}{2q_k}$
of the midpoint of $I$
(recall that $\O_{x_0}$ is dense).  
By the observation in Section \ref{rwonline}, the sequence 
\[\rho(m^*), \rho(f_{\theta_j}(m^*)),\rho(f_{\theta_{j-1}} \circ 
f_{\theta_j}(m^*)), \dots , \rho(f_{\theta_1}\circ \dots \circ 
f_{\theta_j}(m^*))\]  
is a simple random walk on $\Z$.  

If this walk does not stay within $2q_k$ units of $\rho(m^*)$ on the line, 
then we can apply Lemma 
\ref{l:farsmall} to get a $j_0$, $1 \leq j_0 \leq j$ such that 
$f_{\theta_{j_0}}\circ \cdots \circ f_{\theta_{j}}(m^*)<
\frac{3}{2q_k}$. 
Then Lemma \ref{shrinkage} implies that $||f_{\theta_{j_0}}\circ \cdots
\circ f_{\theta_{j}}(I)|| \leq \frac{1}{2}
\length{I} + \delta +\frac{3}{2q_k} \leq 
\frac{1}{2} \length{I} + \frac{2}{q_k} $. 
Since $f_{\theta}$ never
increases the length of an interval, (\ref{intervalshrink}) holds.

With $c_0$ and $N_0$ as in Lemma \ref{l:randwalk}, we therefore have 
$\mathsf{p} 
\geq 1- e^{-2c_0q_k}$ (provided that $2q_k > N_0$).  Thus for 
$n\geq 8q_k^3\log_2 q_k \geq tj$ we have that (\ref{conclusion}) holds
with probability at least $(1-e^{-2c_0q_k})^{\log_2q_k}$.  By the 
binomial theorem, this last expression is the same as
$$1-\frac{\log_2q_k}{e^{2c_0q_k}}+O(\frac{\log^2_2q_k}{e^{4c_0q_k}})$$ 
which, for large $q_k$, is at least $1-c'e^{-2c_0q_k}\log_2q_k>
1-e^{-cq_k}$.  This is the desired result.
\epf

As an application of this theorem, let $\a=1/\sqrt 2$.  Here (see
Section \ref{contfrac}) we have $a_0=0, a_1=1$, and 
$a_i=2$ for $i\geq 2$.  The recurrence relation for $q_n$ 
implies that $$q_n=\frac12((1+\sqrt2)^n+(1-\sqrt2)^n).$$
The theorem implies in this case that there are 
constants $c$ and $C$ such that for $n\geq C\e^{-3}\log(1/\e)$,
and for all $x\in[0,1]$, the $n$th backward iterate is
within $\e$ of its limit with probability at least $1-e^{-c/\e}$.
The same statement holds for any $\a$ for which the set of
$a_i$'s (see (\ref{definea_i})) is bounded; in particular 
this holds for any quadratic irrational $\a$.

\section{Acknowledgements}
We thank Persi Diaconis for suggesting this 
problem and for numerous helpful conversations.  We also 
thank Julie Landau, without whose support and net 
game this paper would not have been written.

\end{document}